\documentclass[reqno,11pt,oneside]{amsart}
\usepackage{amssymb,latexsym}

\newtheorem{theorem}{Theorem}

\newtheorem{lemma}[theorem]{Lemma}
\numberwithin{equation}{section}
\numberwithin{theorem}{section}

\DeclareMathOperator{\sech}{sech}
\newcommand{\Z}{\mathbf{Z}}

\newcommand{\Q}{\mathbf{Q}}
\newcommand{\R}{R}
\newcommand{\ax}{\alpha}
\newcommand{\bx}{\beta}
\newcommand{\n}{\{1,2,\ldots,n\}}
\newcommand{\half}{\tfrac12}

\DeclareMathOperator{\eval}{eval}

\def\e#1#2{e^{(#1)}_{#2}}
\def\f#1#2{f^{(#1)}_{#2}}
\def\r#1#2{(#1)_{#2}}
\def\dpow#1#2{\frac{#1^{#2}}{#2!}}
\def\dpowp#1#2{\frac{#1^{#2}}{(#2)!}}
\def\sumz#1{\sum_{#1=0}^\infty}
\def\egf#1{\sumz n #1 \dpow xn}
\def\egft#1{\sumz n #1 x^n/n!}
\def\ceil#1{\lceil#1\rceil}
\def\qbinom#1#2{\genfrac{[}{]}{0pt}{}{#1}{#2}}

\begin{document}
\title[Applications of the Classical Umbral Calculus]{Applications of the Classical Umbral 
Calculus}
\author{Ira M. Gessel}
\address{Department of Mathematics\\
   Brandeis University\\
   Waltham, MA 02454-9110}
\email{gessel@brandeis.edu}

\date{November 22, 2001} 
\begin{abstract}
We describe  applications of the classical umbral calculus to bilinear generating
functions for polynomial sequences, identities for Bernoulli and related  numbers, and Kummer
congruences. 
\end{abstract}

\dedicatory{Dedicated to the memory of Gian-Carlo Rota}

\thanks{This research was supported by NSF grant DMS-9972648.}
\maketitle

\section{Introduction}  In the nineteenth century, Blissard developed a
notation for manipulating sums involving binomial coefficients by 
expanding polynomials and then  replacing
exponents with subscripts. For example, the expression $(a+1)^n$ would
represent the sum $\sum_{i=0}^n \binom ni a_i$. Blissard's notation has been known variously as
Lucas's method, the symbolic method (or symbolic notation), and  the umbral calculus.
We shall use Rota and
Taylor's term ``classical umbral calculus" \cite{classical} to distinguish it from the more 
elaborate
mathematical edifice that the term ``umbral calculus"  has come to encompass [\ref{roman}, 
\ref{RR}, \ref{RKO}].

The goal of this article is to show, by  numerous examples, how the classical umbral calculus
can be used  to prove interesting
formulas not as easily proved by other methods. Our applications are in three general areas: 
bilinear generating functions, identities for Bernoulli numbers and their relatives,  and
congruences for sequences such as Euler and Bell numbers.

The classical umbral calculus is intimately connected with exponential generating functions; 
thus
$a^n=a_n$ is equivalent to 
\[e^{ax}=\egf {a_n},\] and multiplication of exponential generating functions may be expressed
compactly in umbral notation:
\[\biggl(\egf{a_n}\biggr)\biggl(\egf{b_n}\biggr)=\egf{c_n}\] is equivalent to $(a+b)^n=c^n$. 

When I first encountered umbral notation it seemed to me that this was all there was to
it; it was simply a  notation  for dealing with exponential generating
functions, or to put it bluntly, it was a method for  avoiding the use of exponential
generating functions when they really ought to be used.
The  point of this paper is that my first impression was wrong: none
of the  results proved here (with the exception of Theorem \ref{T:b3}, and perhaps a
few other results in  section
\ref{s:bernoulli}) can be easily proved by straightforward manipulation of exponential 
generating functions. The sequences that we consider here are defined by exponential
generating functions,  and their most fundamental properties can be proved in a
straightforward way  using these  exponential generating functions. What is surprising is
that these sequences satisfy additional relations  whose proofs require other methods. The
classical umbral calculus is a powerful but specialized tool that can be used to prove
these more esoteric formulas. The derangement numbers, for example, have the well-known exponential
generating function
$\sumz n {D_n}x^n/n!=e^{-x}/(1-x)$ from which their basic properties can be derived; 
umbral calculus  gives us the more interesting but considerably more recondite formula
$\sumz n {D_n^2}x^n/n!=e^{x}\sumz k k!\, x^k/{(1+x)^{2k+2}}$.

We begin in the next section with a description of the classical umbral calculus, following Rota
\cite{partitions} and Rota and Taylor  [\ref{umbral}, \ref{classical}], and point out some of
the minor ways in which we differ from their approach. Next, in sections \ref{s:charlier}
through \ref{s:cz}, we consider bilinear and related generating functions for Charlier
and   Hermite polynomials, and some variations. In section \ref{s:rs} we
derive a bilinear generating function for the Rogers-Szeg\H o polynomials, which are related
to $q$-Hermite polynomials.
In section
\ref{s:bernoulli} we apply the umbral calculus to identities for Bernoulli and related numbers.
Sections \ref{s:kummercong} through \ref{s:bell} deal with  Kummer congruences and with 
analogous congruences for Bell
numbers.

The next section contains a formal description of the classical
umbral calculus as used in this paper. The reader who
is not interested in these technicalities may wish to go directly to section
\ref{s:charlier}.

\section{The classical umbral calculus}

Most users of Blissard's symbolic notation  have viewed it as simply a
notational convenience,  requiring no formal justification.  Thus
Guinand \cite{guinand}, in explaining the interpretation of umbral
symbols, writes: ``In general,  any step in manipulation is valid if
and only if it remains valid when interpreted in non-umbral form."
However, in 1940 E. T. Bell \cite{bell} attempted to give an
axiomatic foundation to the umbral calculus. To the modern reader,
Bell's approach seems ill-conceived, if not completely
incomprehensible.  A much more successful explanation was given by
G.-C. Rota in 1964 \cite{partitions}: When we interpret $(a+1)^n$ as
$\sum_{i=0}^n \binom ni a_i$, we are applying the  linear functional
on  the algebra of polynomials in $a$ that takes $a^i$ to~$a_i$. In
retrospect, Rota's idea seems almost obvious, but we must remember that
in Bell's day  the concept of a linear functional was not the familiar
notion that it is in ours. The seemingly mysterious ``umbral variable"
$a$ is just an ordinary variable; it is in the invisible,  but
otherwise unremarkable, linear functional that  the meaning of the
umbral calculus resides. The ``feeling of witchcraft" that Rota and
Taylor \cite{classical} observe hovering about the umbral calculus  comes from
the attribution to umbrae of properties that  really belong to these
linear functionals. As in stage illusion, misdirection is essential
to the magic.

Rota and Taylor's recent works [\ref{umbral}, \ref{classical}]
expanded on Rota's original insight and introduced new concepts
that help to resolve some of the ambiguities that may arise in 
applications of the traditional notation.  However, I shall use the
traditional notation in this paper.
What follows is a short formal description of the classical
umbral calculus as used here, based on Rota and Taylor's
formulation, but with some modifications.

In the simplest applications of the classical umbral calculus, we work
in the ring of polynomials in one variable, e.g., $\R[a]$, where $R$
is a ring of ``scalars" ($\R$ is often a ring of polynomials or formal
power series containing the rationals), and we have a linear
functional $\eval:\R[a]\to \R$. (This notation was introduced by Rota
and Taylor \cite{umbral}.)  The variable
$a$ is called an umbral variable or {\it umbra\/}. There is nothing
special about it other than the fact that the linear functional
$\eval$ is defined on $\R[a]$. 
We will often use the
same letter for the umbra and the sequence; thus we would write $a_n$ for
$\eval(a^n)$.
It is traditional, and
convenient, to omit $\eval$ and to write $a^n=a_n$ instead of
$\eval(a^n)=a_n$. However when following this convention, we must make clear where
$\eval$ is to be applied. The rule that we shall follow in this paper is that $\eval$
should be applied to any term in an equation that contains a variable
that has been declared to be umbral.  It should be emphasized that
this is a syntactic, not mathematical rule, so the formula $a^n=n$ is
to be interpreted as $\eval(a^n)=n$ for all $n$, even though for
$n=0$, $a$ does not ``appear" on the left side.  One important
difference between our approach and that of Rota and Taylor
[\ref{umbral}, \ref{classical}] is that they require that
$\eval(1)=1$, but we do not, and in sections \ref{s:bernoulli} and
\ref{s:eulercong} we shall see several examples where
$\eval(1)=0$. This involves some notational subtleties
discussed below; nevertheless, there is no reason why a linear functional
on polynomials cannot take 1 to 0, and there are are interesting
applications where this happens.

We shall often have occasion to deal with several umbrae together. It should be
pointed out   that although we use the symbol $\eval$ for whatever linear
functional is under discussion, there are really  many different such
functionals. When we  write $a^n=a_n$ and $b^n=b_n$ we are really talking about 
two different linear functionals,
$\eval_1:\R[a]\to\R$ and $\eval_2:\R[b]\to\R$, where $\eval_1(a^n)=a_n$ and
$\eval_2(b^n)=b_n$.  The meaning of $\eval(a^mb^n)$ might be determined by a
completely different linear functional on 
$\R[a,b]$, but traditionally one takes the linear functional $\eval_3$ defined by
$\eval_3(a^m b^n)=\eval_1(a^m)\eval_2(b^n)$. In this case, we say that the umbrae
$a$ and $b$ are {\it  independent\/} (even though we are really dealing with a
property of the linear functional $\eval_3$ rather  than a property of the
variables $a$ and $b$). In fact, applications of umbrae that are not independent in this sense
are uncommon and do not seem to be have been considered before, and we shall
assume that our umbrae are independent  except
where we  explicitly state otherwise.  Nevertheless we give an example
in section \ref{s:cz} of an application of umbrae that are not independent.
 
Eschewing the requirement that $\eval(1)=1$ entails an additional interpretative issue that must 
be
mentioned. We cannot assume that there is a ``universal" evaluation functional that applies to 
every term
in a formula; instead we may need a different functional for each term, corresponding to the 
variables that
appear in that term. In section \ref{s:eulercong}, for example, we have the formula
\[F^n=2A^n-(4B+C)^n,\] involving the umbrae $F$, $A$, $B$, and $C$,
which must be interpreted as
\[\eval_1(F^n)=\eval_2(2A^n) -\eval_3\bigl((4B+C)^n\bigr),\]
where $\eval_1$ is defined on $\Q[F]$, $\eval_2$ is defined on $\Q[A]$, and $\eval_3$ is defined 
on 
$\Q[B,C]$. Although the rule may seem unnatural when stated this way, in practice the 
interpretation is
exactly what one would expect.

We will often find it useful to work with power series, rather than
polynomials, in our umbrae. However, if $f(u)$ is an arbitrary formal
power series in $u$ and $a$ is an umbra then $\eval (f(a))$ does not make sense. Let us
suppose that $\R$ is a ring of formal power series in variables
$x,y,z,\dots$.  Then we call a formal power series $f(u)\in \R[[u]]$
{\it admissible\/} if for every monomial $x^iy^jz^k\cdots$ in $\R$, the
coefficient of $x^iy^jz^k\cdots$ in $f(u)$ is a polynomial
in $u$. Then if $f(u)=\sum_i f_i u^i$ is admissible, we define
$\eval(f(a))$ to be $\sum_i f_i \eval(a^i)$; admissibility of $f$
ensures that this sum is well defined as an element of $\R$. More
generally, we may  define admissibility similarly for a formal power
series in any finite set of variables with coefficients involving
other variables.

\section{Charlier polynomials}
\label{s:charlier}

In the next three sections we apply the classical umbral calculus to find bilinear
generating  functions. More specifically, we find explicit expressions for  generating
functions of the form $\sum_n {a_n b_n}x^n/n!$, where there are simple expressions  for
the generating functions $\sum_n a_n x^n/n!$ and $\sum_n b_n x^n/n!$. Although it is not
obvious  \emph{a priori} that such explicit expressions exist, they do, and they have
important applications in the  theory of orthogonal  polynomials (see, e.g., Askey
\cite{askey}). The method that we use can be translated  into a traditional analytic
computation, since in all cases that we consider in these three sections, $\eval$ can be
represented by a definite integral (though in some cases the radius of convergence of the
series is 0).  For example,  in this section we consider the umbra 
$A$ evaluated by $\eval(A^n)=\ax(\ax+1)\cdots (\ax+n-1)$. We could define $\eval$
analytically  by
\[\eval\bigl(f(A)\bigr)=\frac1{\Gamma(\ax)}\int_0^\infty f(x)x^{\ax -1}e^{-x}\, dx\]
and do all our calculations with integrals. In fact this idea has been used, in a significantly 
more
sophisticated setting, by Ismail and Stanton  [\ref{is1}, \ref{is2}, \ref{is3}] to obtain 
bilinear
generating functions much more complicated than those we deal with here.

 The rising factorial $\r\alpha n$ is defined to be
$\ax(\ax+1)\cdots(\ax+n-1)$. The Charlier polynomials $c_n(x;a)$ are  defined by
\[
c_n(x;a)=\sum_{k=0}^n\binom nk \r{-x}k\, a^{-k}
\]
(see, for example, Askey [\ref{askey}, p.~14]), but 
it is more convenient to work with
differently normalized versions of these polynomials, which we define as 
\[C_n(u,\ax)=u^n c_n(-\ax;u)=\sum_{i=0}^n\binom n i\r\ax{i} u^{n-i}.\]

Let us define the umbra $A$ by $A^n=\r\ax n$. Then
\begin{equation}C_n(u,\ax)=(A+u)^n.\label{eq:C}\end{equation}
 Now 
\begin{equation}e^{Ax}=\sumz n A^n \dpow xn=\sumz n \r\ax n\dpow xn=(1-x)^{-\ax},
\label{eq:0}
\end{equation}
by the binomial theorem. So 
\begin{equation}\sumz n C_n(u,\ax)\dpow xn=e^{(A+u)x}=e^{ux}e^{Ax}=\frac{e^{ux}}{(1-x)^{\ax}}.
\label{eq:1}
\end{equation}

Our goal in this section is to prove the bilinear generating function for
the Charlier polynomials,
\[\sumz n C_n(u,\ax)C_n(v,\bx)\dpow xn=e^{uvx}\sumz k 
\frac{\r \ax k}{(1-vx)^{k+\ax}}
\frac{\r \bx k}{(1-ux)^{k+\bx}}\dpow xk.
\]
To do this we first prove some properties of the umbra $A$.
\begin{lemma}\label{L:1}
For any admissible formal power series $f$, 
\[e^{Ay}f(A)=\frac1{(1-y)^\ax}f\left(\frac A{1-y}\right).\]
\end{lemma}

\begin{proof}
First we prove the lemma for the case $f(z)=e^{zw}$. We have
\begin{align}
e^{Ay}e^{Aw}&=e^{A(y+w)}=\frac1{(1-y-w)^\ax}\notag\\ 
&=\frac1{(1-y)^\ax}
\frac1{\left(\displaystyle 1-\frac w{1-y}\right)^\ax}\notag\\
&=\frac1{(1-y)^\ax}\exp\left(\frac A{1-y}w\right).
\label{eq:A}\end{align}
by \eqref{eq:0}. Equating coefficients of $w^k\!/k!$ shows that the lemma is true for 
$f(z)=z^k$. The general case then follows by linearity.

Alternatively, we could have introduced an umbra $F$ with $e^{Fz}=f(z)$ and
replaced $w$ with  $F$ in \eqref{eq:A}.
\end{proof}

As a first application of Lemma \ref{L:1}, we prove the following
little-known result.

\begin{theorem}\label{T:double}
\[\sumz m C_{2m}(u,\ax)\dpow xm=e^{u^2x}
\sumz k \frac{\r \ax{2k}}{(1-2ux)^{2k+\ax}}\dpow xk.\]
\end{theorem}

\begin{proof}
We have
\begin{align*}
\sumz m C_{2m}(u,\ax)\dpow xm 
  &= \sumz m(A+u)^{2m}\dpow xm=e^{(A+u)^2x}\\
  &= e^{(A^2+2Au+u^2)x}=e^{u^2x}e^{2Aux}e^{A^2x}\\
  &=\frac{e^{u^2x}}{(1-2ux)^\ax}    
   \exp\left[\left(\frac A{1-2ux}\right)^2 x\right]
   \text{\quad by Lemma \ref{L:1}}\\
  &=e^{u^2 x}\sumz k \frac{\r \ax{2k}}{(1-2ux)^{2k+\ax}}\dpow xk.
\end{align*}
\end{proof}

By a similar computation we can prove a generalization given by the next
theorem. We leave the details to the reader.

\begin{theorem}
\[\sumz{m,n}C_{2m+n}(u,\ax) \dpow xm \dpow yn=
e^{u^2x+uy}\sumz k \frac{\r\ax{2k}}{(1-2ux-y)^{2k+\ax}}\dpow xk.\qed\]
\end{theorem}

Next we prove the bilinear generating function for Charlier polynomials. An equivalent formula 
can be found in Askey
[\ref{askey}, p.~16, equation (2.47)] with a minor error; $a$ and $b$ must be switched on one 
side of the
formula as given there for it to be correct. A combinatorial proof of our Theorem 
\ref{T:bilC} has been given
by Jayawant \cite{jayawant}, who also proved a multilinear generalization.
\begin{theorem}
\label{T:bilC}
\[\sumz n C_n(u,\ax)C_n(v,\bx)\dpow xn=e^{uvx}\sumz k 
\frac{\r \ax k}{(1-vx)^{k+\ax}}
\frac{\r \bx k}{(1-ux)^{k+\bx}}\dpow xk.
\]
\end{theorem}
\begin{proof} Let $A$ and $B$ be independent umbrae with 
$A^n=\r \ax n$ and $B^n=\r \bx n$.  Then there is an analogue of Lemma
\ref{L:1} with $B$ replacing $A$ and $\bx$ replacing $\ax$.

We have
{\allowdisplaybreaks
\begin{align*}
\sumz n C_n(u,\ax)C_n(v,\bx)\dpow xn&=e^{(A+u)(B+v)x}\\
    &=e^{uvx}e^{Avx}e^{(Bu+AB)x}\\
    &=e^{uvx}\frac1{(1-vx)^\ax}\exp\left(Bux +\frac A{1-vx}Bx\right) 
           \text{\quad by Lemma \ref{L:1}}\\
    &=\frac{e^{uvx}}{(1-vx)^\ax}e^{Bux}
       \exp\left(\frac A{1-vx}Bx\right)\\
    &=\frac{e^{uvx}}{(1-vx)^\ax}\cdot
       \frac1{(1-ux)^\bx}
       \exp\left(\frac A{1-vx}\cdot\frac B{1-ux}x\right)
        \text{\quad by Lemma \ref{L:1}}\\
    &=e^{uvx}\sumz k 
      \frac{\r \ax k}{(1-vx)^{k+\ax}}
      \frac{\r \bx k}{(1-ux)^{k+\bx}}\dpow xk.
\end{align*}
}

\end{proof}

The polynomials $C_n(u,\ax)$ have a simple interpretation in terms of
permutation enumeration: the coefficient of  $\ax^i u^j$ in 
$C_n(u-\ax,\ax)$ is the number of permutations of $\n$ with 
$j$ fixed points and $i$ cycles of length at least 2. This follows
easily from the exponential generating function
\[e^{ux}\left(\frac{e^{-x}}{1-x}\right)^\ax=\sumz n C_n(u-\ax,\ax)\dpow
xn.\] (See, for example, Stanley [\ref{stanley}, chapter 5].) In particular, $C_n(-1,1)$ is
the derangement number $D_n$, the number of permutations of $\n$ with no
fixed points, and Theorems \ref{T:double} and \ref{T:bilC} give the
formulas
\[\sumz m D_{2m}\dpow xm=e^{x}\sumz k\frac{(2k)!}
{(1+2x)^{2k+1}}\dpow xk\]
and 
\[\sumz n D_n^2\dpow xn=e^{ x}\sumz k \frac{k!}{(1+x)^{2k+2}}x^k.\]

Theorem \ref{T:bilC} can be generalized to a formula involving 3-line
Latin rectangles. See \cite{3line} for a combinatorial proof that also
uses umbral methods. A more general result  was given using the same technique by Zeng 
\cite{zeng1}, and using very
different techniques by Andrews, Goulden, and Jackson \cite{agj}.

\section{Hermite polynomials}
\label{s:hermite}
We now prove some similar formulas for Hermite polynomials. Perhaps
surprisingly, the proofs are a little harder than those for Charlier
polynomials. We first define the umbra $M$ by
\begin{equation}e^{Mx}=e^{-x^2},\label{eq:M}\end{equation} so that
\[M^n=\begin{cases}\displaystyle(-1)^k\frac{(2k)!}{ k!},&
\text{if
$n=2k$}\\ 0,& \text{if $n$ is odd.}\end{cases}
\]
(The  reason for the minus sign in this definition  is so that we can obtain
formulas for the Hermite polynomials in their usual normalization.) There are two basic
simplification formulas for
$M$:

\begin{lemma}
\label{L:M}\ \par
\begin{enumerate}
\item[(i)] $\displaystyle e^{M^2x}=\frac 1{\sqrt{1+4x}}.$
\openup 2\jot
\item[(ii)]
For any admissible formal power series $f$, we have
\[e^{My}f(M)=e^{-y^2}f(M-2y).\]
\end{enumerate}\end{lemma}

\begin{proof}
For (i), we have
\[e^{M^2x}=\sumz k (-1)^k\frac{(2k)!}{ k!}\dpow xk=\frac1{\sqrt{1+4x}}.\]
For (ii), as in the proof of Lemma \ref{L:1}, it is sufficient to prove
that the formula holds for $f(z)=e^{zw}$. In this case we have
\[e^{My}e^{Mw}=e^{M(y+w)}=e^{-y^2-2yw-w^2}=e^{-y^2}e^{-2yw}e^{Mw}
  =e^{-y^2}e^{(M-2y)w}.\]
\end{proof}

\begin{lemma}
\label{L:Mquad}
\[e^{Mx+M^2y}=\frac{e^{-x^2/(1+4y)}}{\sqrt{1+4y}}.\]
\end{lemma}
Although Lemma \ref{L:Mquad} can be proved directly by showing that both sides are equal to 
\[\sum_{i,j}(-1)^{j+k}\frac{(2j+2k)!}{(j+k)!}\dpowp x{2j}\dpow yk,\] 
we  give instead two proofs that use Lemma \ref{L:M}
\begin{proof}[First proof]
If we try to apply Lemma \ref{L:M} directly, we find that the linear term in $M$ does not 
disappear, so we need to use a
slightly less direct approach. We write $e^{Mx+M^2y}$ as $e^{M(x+z)}e^{M^2y-Mz}$, where $z$ will 
be chosen later. Now
applying Lemma \ref{L:M} gives
\begin{align*} e^{M(x+z)}e^{M^2y-Mz}&=e^{-(x+z)^2}e^{(M-2x-2z)^2y-(M-2x-2z)z}\\
&=e^{-(x+z)^2}e^{M^2y-4M(x+z)y+4(x+z)^2y-Mz+(2x+2z)z}.\end{align*}
We now choose $z$ so as to eliminate the linear term in $M$ on the right; i.e., we want
$-4(x+z)y-z=0$. So we take $z=-4xy/(1+4y)$, and on simplifying we obtain 
$e^{Mx+M^2y}=e^{-x^2/(1+4y)+M^2y}$. Then applying Lemma \ref{L:M} (i) gives the desired
result.
\end{proof}

\begin{proof}[Second proof]
Let us fix $y$ and set $g(x)=e^{Mx+M^2y}$.
Applying Lemma \ref{L:M} directly gives
\begin{align*}
g(x)&=e^{Mx+M^2y}
  =e^{-x^2}e^{(M-2x)^2 y}\\
  &=e^{-x^2+4x^2y} e^{-4Mxy+M^2y}
  =e^{-x^2(1-4y)} g(-4xy).
\end{align*}
Iterating and taking a limit yields
\begin{align*}
g(x)&=e^{-x^2(1-4y) -4^2x^2y^2(1-4y)-\cdots}
    =e^{-x^2(1-4y+4^2y^2-4^3y^3+\cdots)}g(0)\\
    &=e^{-{x^2}/(1+4y)}g(0)=e^{-x^2/(1+4y)}/\sqrt{1+4y}
\end{align*}
by Lemma \ref{L:M} (i).
\end{proof}

Now we define the Hermite polynomials $H_n(u)$ by the generating function
\begin{equation}\sumz n H_n(u)\dpow xn=e^{2ux-x^2}=e^{(2u+M)x}
\label{e:hermite}
\end{equation}
so that $H_n(u)=(2u+M)^n$.

First we prove a well-known analogue of Theorem \ref{T:double}, a special case of a result of Doetsch \cite[equation
(10)]{doetsch}.
  
\begin{theorem}\label{t:doetsch}
\[\sumz n H_{2n}(u)\dpow xn
=\frac1{\sqrt{1+4x}}\exp\left(\frac{4u^2 x}{1+4x}\right).
\]
\end{theorem}

\begin{proof}
We have
\begin{align*}
\sumz n H_{2n}(u)\dpow xn
 &=e^{(2u+M)^2 x}=e^{4u^2x}e^{4Mux+M^2x}\\
 &=\frac1{\sqrt{1+4x}}e^{4u^2x}\exp\left(\frac{-16u^2 x^2}{1+4x}\right)
 \text{\quad by Lemma \ref{L:Mquad}}\\
 &=\frac1{\sqrt{1+4x}}\exp\left(\frac{4u^2 x}{1+4x}\right).
\end{align*}
\end{proof}

By the same reasoning we can prove the following generalization of Theorem 
\ref{t:doetsch}.

\begin{theorem}
\label{t:doetsch2}
\[\sumz {m,n} H_{2m+n}(u) \dpow xm\dpow yn=
\frac1{\sqrt{1+4x}}\exp\left(\frac{4u^2x+2uy-y^2}{1+4x}\right).\qed\] 
\end{theorem}

Equating coefficients of $y^n/n!$ in both sides of Theorem \ref{t:doetsch2}, and using
\eqref{e:hermite} yields
 \[\sumz m H_{2m+n}(u)\dpow xm
= (1+4x)^{-(n+1)/2}H_n\left(\frac u {\sqrt{1+4x}}\right)
\exp\left(\frac{4u^2x}{1+4x}\right),
\]
which is the general form of Doetsch's result \cite{doetsch}.

We state without proof a ``triple" version of Theorem \ref{t:doetsch} that can be proved  by the 
same
technique. See Jayawant \cite{jayawant}, where umbral and combinatorial proofs are given. 

\begin{theorem}
\[\egf{H_{3n}(u)}=\frac{e^{8v^3 x+144 v^4 x^2}}{(1+48ux)^{1/4}}\sumz 
n\frac{(-1)^n(6n)!}{(3n)!\,(1+48ux)^{3n/2}}
\dpowp x{2n},\]
where 
$v=\left(\sqrt{1+48ux}-1\right)/(24x)$. \qed
\end{theorem}

%
%
%

\bigskip
Next we prove Mehler's formula, which gives a bilinear generating function for the
Hermite polynomials. An elegant combinatorial proof of this formula has been given by
Foata \cite{foata1}, and generalized to the multilinear case by Foata and Garsia
[\ref{foata2}, \ref{fogar}].

\begin{theorem}
\label{mehler}
\[\egf {H_n(u) H_n(v)}=\frac1{\sqrt{1-4x^2}}\exp\left(4\frac{uvx
-(u^2+v^2)x^2}{1-4x^2}\right).\]
\end{theorem}
\begin{proof}
We use two independent umbrae, $M$ and $N$, with $M$ as before and $N^n=M^n$ for
all~$n$. (In the terminology of Rota and Taylor [\ref{classical}, \ref{umbral}],
$M$ and $N$ are ``exchangeable" umbrae.) Then
\begin{align*}
\egf {H_n(u) H_n(v)} &=e^{(2u+M)(2v+N)x}\\
  &=e^{2u(2v+N)x}e^{M(2v+N)x}\\
  &=e^{2u(2v+N)x}e^{-(2v+N)^2x^2}\text{\quad by \eqref{eq:M}}\\
  &=e^{4vx(u-vx)}e^{2Nx(u-2vx)-N^2x^2}\\
  &=\frac{e^{4vx(u-vx)}}{\sqrt{1-4x^2}}\exp\left(-\frac{4x^2(u-2vx)^2}{1-4x^2}\right)
       \text{\quad by Lemma \ref{L:Mquad}}\\
  &=\frac1{\sqrt{1-4x^2}}\exp\left(4\frac{uvx -(u^2+v^2)x^2}{1-4x^2}\right).
\end{align*}
\end{proof}

\section{Carlitz and Zeilberger's Hermite polynomials}
\label{s:cz}
Next we consider analogues of the Hermite polynomials studied by Carlitz \cite{carlitz3} and 
Zeilberger
\cite{zeilberger}. Carlitz considered the ``Hermite polynomials of two variables"
\[H_{m,n}(u,v)=\sum_{k=0}^{\min(m,n)}\binom mk \binom nk k!\, u^{m-k} v^{n-k},\]
with generating function
\begin{equation*}\sum_{m,n}H_{m,n}(u,v) \dpow xm \dpow yn=e^{ux+vy+xy}
\end{equation*}
and proved the  bilinear generating function
\begin{multline}
\sumz{m,n} H_{m,n}(u_1,v_1)H_{m,n}(u_2,v_2)\dpow xm \dpow yn\\
=(1-xy)^{-1}\exp\left(\frac{u_1u_2x +v_1v_2y+(u_1 v_1+u_2 v_2)xy}{1-xy}\right).
\label{e:carlitz}
\end{multline}
Independently, Zeilberger considered the ``straight Hermite polynomials"
\[H_{m,n}(w)=\sum_{k=0}^{\min(m,n)}\binom mk\binom nk k!\, w^k,\]
with generating function 
\[
  \sum_{m,n}H_{m,n}(w) \dpow xm \dpow yn=e^{x+y+wxy},
\]
and gave a combinatorial proof, similar to Foata's proof of Mehler's formula \cite{foata1}, of 
the bilinear generating
function
\begin{equation}
\sumz{m,n} H_{m,n}(u)H_{m,n}(v)\dpow xm \dpow yn
=(1-uvxy)^{-1}\exp\left(\frac{x+y+(u+v)xy}{1-uvxy}\right).\label{e:zeil}
\end{equation}
It is easy to see that Carlitz's and Zeilberger's polynomials are related by 
$H_{m,n}(u,v)=u^m v^n H_{m,n}(1/uv)$, and that \eqref{e:carlitz} and \eqref{e:zeil} are 
equivalent. We shall prove
\eqref{e:zeil}, since it involves fewer variables. Our proof uses umbrae that are not 
independent.

We define the umbrae $A$ and $B$ by
\[A^m B^n = \delta_{m,n} m!,\]
where $\delta_{m,n}$ is $1$ if $m=n$ and $0$ otherwise. Equivalently, $A$ and $B$ may be defined 
by 
\begin{equation}
e^{Ax+By}=e^{xy}.
\label{e:AB1}
\end{equation}
Then Zeilberger's
straight Hermite polynomials are given by
$H_{m,n}(u)=(1+A)^m(1+Bu)^n$. Two of the basic properties of these umbrae are given in the 
following
lemma.

\begin{lemma}
\label{L:AB1}
\ \par
\begin{enumerate}
\item[(i)] If $f(x,y)$ is an admissible power series then $e^{Ar+Bs}f(A,B)=e^{rs}f(A+r,B+s)$.
\openup 2\jot
\item[(ii)]
$\displaystyle e^{Ax+By +ABz}=\frac1{1-z}e^{xy/(1-z)}.$
\end{enumerate}
\end{lemma}

\begin{proof}
As in the proof of Lemma \ref{L:1}, it is sufficient to prove (i) for the case 
$f(x,y)=e^{xu+yv}$, and for this case we
have
\[e^{Ar+Bs}f(A,B)=e^{Ar+Bs}e^{Au+Bv}=e^{A(r+u)}e^{B(s+v)}=e^{(r+u)(s+v)},\] 
by \eqref{e:AB1},
and 
\begin{align*}
e^{rs}f(A+s,B+r)&=e^{rs}e^{(A+s)u+(B+r)v}=e^{rs+rv+su}e^{Au+Bv}\\
&=e^{rs+rv+su}e^{uv}=e^{(r+u)(s+v)}.
\end{align*}
We can prove (ii) by using (i) to reduce it to the case $x=y=0$, but instead we give a direct 
proof.
We have
\begin{align*}
e^{Ax+By+ABz}&=\sum_{i,j,k}A^{i+k}B^{j+k}\dpow xi \dpow yj\dpow zk\\
  &=\sum_{j,k}(j+k)!\, \frac{(xy)^j}{j!^2}\dpow zk
    =\sum_{j,k}\binom{j+k}j \dpow{(xy)}j z^k\\
  &=\sum_j\frac1{(1-z)^{j+1}}\dpow{(xy)}j
    =\frac1{1-z}e^{xy/(1-z)}.
\end{align*}
\end{proof}
Now we prove Zeilberger's bilinear generating function \eqref{e:zeil}. We introduce two 
independent pairs of umbrae
$A_1,B_1$ and
$A_2,B_2$ such that each pair behaves like $A,B$; in other words,
\[A_1^k B_1^l A_2^m B_2^n=\delta_{k,l}\delta_{m,n}k!\,m!.\]

Then 
\begin{align*}
\sumz{m,n} H_{m,n}(u)&H_{m,n}(v)\dpow xm \dpow yn\\
  &=\sum_{m,n}(1+A_1)^m(1+B_1 u)^n (1+A_2)^m(1+B_2 v)^n \dpow xm \dpow yn\\
  &=e^{(1+A_1)(1+A_2)x+(1+B_1 u)(1+ B_2 v)y}\\
  &=e^{ (1+A_2)x+(1+B_2 v)y} e^{A_1(1+A_2)x+B_1(1+B_2 v)uy}.
\end{align*}
Applying \eqref{e:AB1} with $A_1$ and $B_1$ for $A$ and $B$ yields
\begin{equation*}
e^{(1+A_2)x+(1+B_2 v)y+(1+A_2)(1+B_2 v)uxy} 
  = e^{x+y+uxy} e^{A_2x(1+uy)+B_2vy(1+ux)+A_2B_2uvxy}.
\end{equation*}
Then applying Lemma \ref{L:AB1} (ii) yields \eqref{e:zeil}.\medskip

By similar reasoning, we can prove a generating function identity equivalent to the 
Pfaff-Saalsch\"utz theorem for
hypergeometric series \cite{gs}. In terms of Carlitz's Hermite polynomials of two variables, 
this is the
evaluation of \[\sum_{m,n}H_{m,n+j}(0,1)H_{m+i,n}(0,1) \dpow xm \dpow yn.\]
\begin{theorem}
Let $i$ and $j$ be nonnegative integers. Then 
\begin{equation}\label{e:saal}
\sumz{m,n}\binom{m+i}n\binom{n+j}m x^m y^n=\frac{(1+x)^j(1+y)^i}{(1-xy)^{i+j+1}}.
\end{equation}
\end{theorem}
\begin{proof}
With $A_1$, $B_1$, $A_2$, and $B_2$ as  before, we have
\[A_1^m(1+B_1)^{n+j}A_2^n(1+B_2)^{m+i}=m!\,n!\,\binom{m+i}n\binom{n+j}m,\] 
so the left side of \eqref{e:saal} is equal to 
\begin{multline}
\label{e:saal2}
\qquad\sum_{m,n}A_1^m(1+B_1)^{n+j}A_2^n (1+B_2)^{m+i}\dpow xm\dpow yn\\
=e^{A_1(1+B_2)x+A_2(1+B_1)y}(1+B_2)^i(1+B_1)^j.\qquad
\end{multline}
Multiplying the right side of \eqref{e:saal2} by $u^iv^j/i!\,j!$, and summing on $i$ and $j$, we 
obtain
\[
e^{A_1(1+B_2x)+A_2(1+B_1)y+(1+B_2)u+(1+B_1)v}
  =e^{u+v}e^{A_1(1+B_2)x+B_1(v+A_2y)+A_2y+B_2u}.
\]
Applying \eqref{e:AB1}, with $A_1$ and $B_1$ for $A$ and $B$, gives
\begin{align*}
e^{u+v}e^{(1+B_2)(v+A_2y)x+A_2y+B_2u}
  &=e^{u+v+xv}e^{A_2 (1+x)y+B_2(u+xv)+A_2B_2xy}.
\end{align*}
Applying Lemma \ref{L:AB1} (ii), we obtain
\[
\frac{e^{u+v+xv}}{1-xy}\exp\left(\frac{(1+x)(u+xv)y}{1-xy}\right)
=\frac1{1-xy}\exp\left(\frac{(1+x)v+(1+y)u}{1-xy}\right),
\]
and extracting the coefficient of  $u^iv^j/i!\,j!$ gives the desired result.
\end{proof}

We can also prove analogues of Doetsch's theorem (Theorem \ref{t:doetsch}) for the straight 
Hermite polynomials. We need
the following lemma, which enables us to evaluate the exponential of any quadratic
polynomial in $A$ and $B$.

\begin{lemma}
\label{L:AB2}
\[e^{Av+Bw+A^2x+ABy+B^2z}=\frac1{\sqrt{(1-y)^2-4xz}}\exp\left(\frac{vw(1-y)+v^2z+w^2x}{(1-y)^2-4
xz}\right).\]
 \end{lemma}

\begin{proof}
Since the proof is similar to earlier proofs, we omit some of the details. The case $v=w=0$ is 
easy to prove directly.
For the general case, we write $e^{Av+Bw+A^2x+ABy+B^2z}$ as $e^{Ar+Bs}\cdot 
e^{-Ar-Bs+Av+Bw+A^2x+ABy+B^2z}$ and choose
$r$ and $s$ so that when Lemma \ref{L:AB1} (i) is applied, the linear terms in $A$ and $B$ 
vanish. We find that the right
values for $r$ and $s$ are
\[r=\frac{v(1-y)+2wx}{(1-y)^2-4xz}\text{\quad and \quad}s=\frac{w(1-y)+2vz}{(1-y)^2-4xz},\]
and the result of the substitution is 
\[\exp\left(\frac{vw(1-y)+v^2z+w^2x}{(1-y)^2-4xz}\right) e^{A^2x+ABy+B^2z},\]
which may be evaluated by the case $v=w=0$.
\end{proof}

\begin{theorem}
\label{T:Zdoetsch}
\begin{align*}
\sumz{m,n}H_{2m,n}(u)\dpow xm\dpow yn&=e^{x+y+2uxy+u^2xy^2}\\
\sumz{m,n}H_{2m,2n}(u)\dpow xm\dpow 
yn&=\frac1{\sqrt{1-4u^2xy}}\exp\left(\frac{x+y+4uxy}{1-4u^2xy}\right)\\
\sumz m H_{m,m}(u)\dpow xm&=\frac1{1-ux}\exp\left(\frac x{1-ux}\right)
\end{align*}
\end{theorem}
\begin{proof}
For the first formula, we have 
\[\sumz{m,n}H_{2m,n}(u)\dpow xm\dpow yn=\sumz{m,n}(1+A)^{2m}(1+uB)^n\dpow xm\dpow yn
=e^{(1+A)^2 x +(1+uB)y}.\]
We simplify this with Lemma \ref{L:AB2}. The proofs of the other two formulas are similar. (The 
third formula
is equivalent to a well-known generating function for Laguerre polynomials.)
\end{proof}

By the same reasoning, we can prove a more general formula that includes all three formulas
of  Theorem \ref{T:Zdoetsch} as special cases.

\begin{theorem}
\begin{multline}
\sumz {i,j,k,l,m}H_{i+2k+m,j+2l+m}(u)\dpow vi\dpow wj\dpow xk\dpow yl\dpow zm
=
\frac1{\sqrt{(1-uz)^2-4u^2xy}}\\
\times
\exp\left(\frac
{(1+uw)^2 x +(1+uv)^2 y+4uxy+(1-uz)(v+w+z+uvw)}
{(1-uz)^2-4u^2xy}
\right).\qed
\end{multline}
\end{theorem}

\section{Rogers-Szeg\H o polynomials}
\label{s:rs}
Next we give a proof of a bilinear generating function for the Rogers-Szeg\H o polynomials, 
which are closely related
to $q$-Hermite polynomials. Our proof differs from the other proofs in this paper in that it 
uses a linear functional
on a noncommutative polynomial algebra. A traditional proof of this result can be found in 
Andrews [\ref{andrews}; p.~50, Example 9] which is also a good reference for basic facts about
$q$-series.

In this section we let $(a)_m$ denote the $q$-factorial
\[(a)_m=(1-a)(1-aq)\cdots (1-aq^{m-1}),\] with $(a)_\infty=\lim_{m\to\infty}\r am$ as a
power  series in $q$. In
particular,
\[(q)_m=(1-q)(1-q^2)\cdots(1-q^m).\]  The
$q$-binomial coefficient  $\qbinom nk$ is defined to be $(q)_n/(q)_k(q)_{n-k}$.
The {Rogers-Szeg\H o} polynomials
$R_n(u)$ are defined by
\[R_n(u)=\sum_{k=0}^n \qbinom n k u^k.\]
We will use a $q$-analogue of the exponential function,
\[e(x)=\sumz n \frac{x^n}{\r qn}.\]
We will also need the $q$-binomial theorem
\[\sumz n \frac{\r an}{\r qn}x^n=\frac{\r {ax}\infty}{\r x\infty};\]
the special case $a=0$ gives
\[e(x)=\frac1{(x)_\infty},\] from which it follows
that 
$e(q^jx)=(x)_je(x)$.

If $A$ and $B$ are noncommuting variables
satisfying the commutation relation $BA=qAB$, then it is well known that 
\begin{equation}(A+B)^n=\sum_{k=0}^n \qbinom nk A^k B^{n-k},\label{eq:nonc}\end{equation}
and it follows easily from \eqref{eq:nonc} that
$e\bigl((A+B)x\bigr)=e(Ax)e(Bx)$, where $x$ commutes with $A$ and $B$.
We shall also need the easily-proved fact that $B^jA^i=q^{ij}A^iB^j$.
 
Now let $A,B,C$, and $D$ be noncommuting variables such that $BA=qAB$,
$DC=qCD$, and all other pairs of variables commute. We shall work in the ring of
formal power series in $A$, $B$, $C$, $D$, with our ring of scalars
(which commute with everything) containing variables $u$, $v$, $x$ and
$q$.  We define our evaluation functional by $\eval(A^iB^jC^kD^l)=u^iv^k$.

Since we need to do some of our computations in the ring of formal power series in $A$, $B$, 
$C$, and $D$,  we write out
the applications of $\eval$ explicitly in this proof.
\begin{theorem}
\[\sumz n R_n(u) R_n(v)\frac{x^n}{\r qn}=\frac{(uvx^2)_\infty}{(uvx)_\infty(ux)_\infty
(vx)_\infty(x)_\infty}.\]
\end{theorem}

\begin{proof}
By \eqref{eq:nonc},
\begin{equation}
\label{E:q1}\eval\left(\sumz n(A+B)^n(C+D)^n\frac{x^n}{\r qn}\right)=\sumz n
R_n(u)R_n(v)\frac{x^n}{\r qn}.\end{equation}
Also, we have
\begin{align}
\sumz n(A+B)^n(C+D)^n\frac{x^n}{\r qn}&=e\bigl((A+B)(C+D)x\bigr)\notag \\
&=e\bigl(A(C+D)x\bigr)e\bigl(B(C+D)x\bigr)\notag\\
&=e(ACx)e(ADx)e(BCx)e(BDx).\label{E:q2}
\end{align}
The only variables ``out of order" in this product are the $D$'s and $C$'s in 
$e(ADx)e(BCx)$,
so
\begin{align}\eval\bigl(e(ACx)e(ADx)e(BCx)&e(BDx))\notag\\
  &=\eval\bigl(e(ACx)\bigr)\eval\bigl(e(ADx)e(BCx)\bigr)\eval\bigl(e(BDx)\bigr)\notag\\
  &=e(uvx)\eval\bigl(e(ADx)e(BCx)\bigr)e(x)\notag\\
  &=\frac{\eval\bigl(e(ADx)e(BCx)\bigr)}
            {\r x\infty\r{uvx}\infty}.\label{E:q3}
\end{align}
We have
\begin{equation*}
e(ADx)e(BCx)=\sumz{i,j}\frac{(ADx)^i}{\r qi}\frac{(BCx)^j}{\r qj} 
=\sumz {i,j}\frac{A^iB^jC^jD^i q^{ij}x^{i+j}}{\r qi\r qj},
\end{equation*}
so 
\begin{align}
\eval\bigl(e(ADx)e(BCx)\bigr)&=\sumz {i,j}\frac{u^iv^jq^{ij}x^{i+j}}{\r qi\r qj}
  =\sumz i\frac{(ux)^i}{\r qi}\sumz j\frac{(vxq^i)^j}{\r qj}\notag\\
  &=\sumz i\frac{(ux)^i}{\r qi \r {vxq^i}\infty}
  =\frac1{\r {vx}\infty}\sumz i\frac{\r{vx}i}{\r qi}(ux)^i\notag\\
  &=\frac1{\r {vx}\infty}\frac{\r{uvx^2}\infty}{\r {ux}\infty}.\label{E:q4}
\end{align}

The theorem then follows from \eqref{E:q1}, \eqref{E:q2}, \eqref{E:q3}, and \eqref{E:q4}.
\end{proof}

It is worth pointing out that although our proof uses noncommuting variables, it does not yield 
a
noncommutative generalization of the result, since the last application of $\eval$ is necessary 
for the final
simplification.

\section{Bernoulli numbers}
\label{s:bernoulli}
The Bernoulli numbers $B_n$  are defined by the exponential generating function

\begin{equation}
B(x)=\egf{B_n}=\frac x{e^x-1}.\label{e:bgf}
\end{equation}
Since \eqref{e:bgf} implies that 
$e^xB(x)=x+B(x)$, the Bernoulli umbra $B$ defined by $B^n=B_n$ satisfies
\begin{equation}
(B+1)^n=B^n+\delta_{n-1},
\label{e:b1}
\end{equation} where
$\delta_m$ is 1 if $m=0$ and is 0 otherwise. From \eqref{e:b1} it follows by linearity that for 
any admissible formal power
series $f$, 
\begin{equation}
f(B+1)=f(B)+f'(0).
\label{e:b2}
\end{equation}
Formula \eqref{e:b2} may be iterated to yield 
\begin{equation}
f(B+k)=f(B)+f'(0) +f'(1)+\cdots+f'(k-1)
\label{e:b3}
\end{equation}
for any nonnegative integer $k$.

There are three other important basic identities for the Bernoulli
umbra. Although the most straightforward proofs use exponential
generating functions, the umbral proofs are interesting and are therefore
included here. Very different umbral proofs of these identities have
been given by Rota and Taylor [\ref{classical}, Theorem 4.2 and
Proposition 8.3].

\begin{theorem}
\label{T:b3}\ \par
\begin{enumerate}
\item[(i)] $(B+1)^n=(-B)^n.$ 
\openup 2\jot
\item[(ii)] $(-B)^n=B^n $ for $n\ne 1$, with $B_{1}=-\tfrac12$. Thus $B_n=0$ when $n$ is odd and
greater than~1.
\item[(iii)]
For any positive integer $k$,
\[\displaystyle kB^n=(kB)^n+(kB+1)^n+\cdots + (kB+k-1)^n.\]
\end{enumerate}
\end{theorem}

\begin{proof} We prove ``linearized" versions of these formulas: for any polynomial $f$,
we  have
\begin{align}
f(B+1)&=f(-B)\label{e:bx1}\\
f(-B)&=f(B)+f'(0)\label{e:bx2}\\
kf(B)&=f(kB)+f(kB+1)+\cdots+f(kB+k-1)\label{e:bx3}
\end{align}
First note that \eqref{e:bx2} follows immediately from \eqref{e:bx1}
and \eqref{e:b2}.  We prove \eqref{e:bx1} and \eqref{e:bx3} by
choosing polynomials $f(x)$, one of each possible degree, for which
the formula to be proved is an easy consequence of \eqref{e:b2}.

For \eqref{e:bx1}, we take
$f(x)=x^n-(x-1)^n$, where
$n\ge 1$.
Then \[f(B+1)=(B+1)^n-B^n=\delta_{n-1}\text{\quad by \eqref{e:b1}},\] and 
since $f(-x)=(-1)^{n-1}f(x+1)$, we have
\[f(-B)=(-1)^{n-1}f(B+1)=(-1)^{n-1}\delta_{n-1}=\delta_{n-1}=f(B+1).\]

For \eqref{e:bx3}, we take $f(x)=(x+1)^n-x^n$, where $n\ge1$. Then 
$f(B)=\delta_{n-1}$ and
\begin{align*}
\sum_{i=0}^{k-1}f(kB+i)&=\sum_{i=0}^{k-1}(kB+i+1)^n-\sum_{i=0}^{k-1}(kB+1)^n\\
  &=(kB+k)^n -(kB)^n = k^n\bigl((B+1)^n-B^n\bigr)\\
  &=k^n\delta_{n-1}=k\delta_{n-1}=kf(B).
\end{align*}
\end{proof}

For later use, we note two consequences of Theorem \ref{T:b3}. First, 
combining \eqref{e:b3} and \eqref{e:bx2} gives
\begin{equation}
\label{e:b4}
f(B+k)-f(-B)=\sum_{i=1}^{k-1}f'(i).
\end{equation}
Second, suppose that $f(u)$ is a polynomial satisfying $f(u+1)=f(-u)$.
Then we have
\begin{align}
f(B)&=\half\bigl(f(2B)+f(2B+1)\bigr)\quad\text{by \eqref{e:bx3}} \notag\\
&=\half\bigl(f(2B)+f(-2B)\bigr)\notag\\
&=f(2B)+f'(0)\quad\text{by \eqref{e:bx2}}.\label{e:double}
\end{align}

\def\B{\tilde B}
Next, we discuss an identity of Kaneko \cite{kaneko}, who set $\B_n=(n+1)B_n$ and gave
the identity 
\begin{equation}
\label{e:Kaneko}
\sum_{i=0}^{n+1}\binom{n+1}i \B_{n+i}=0,
\end{equation}
noting that it (together with the fact that $B_{2j+1}=0$ for $j>0$) allows the computation
of $B_{2n}$ from only half of the preceding Bernoulli numbers. Kaneko's proof is
complicated, though his paper also contains a short proof by D. Zagier. We shall
show that Kaneko's identity is a consequence of the following nearly trivial result.

\begin{lemma}
\label{L:Kaneko}
For any nonnegative integers $m$ and $n$,
\begin{equation*}
\sum_{i=0}^m\binom mi B_{n+i}=(-1)^{m+n}\sum_{j=0}^n\binom nj B_{m+j}.
\end{equation*}
\end{lemma}
\begin{proof}
Take $f(x)=x^m(x-1)^n$ in \eqref{e:bx1}.
\end{proof}

The key to Kaneko's identity is the observation that
\begin{equation}
\label{e:Kkey}
\binom {n+1}i\B_{n+i}=(n+1)\left[\binom{n+1}i +\binom n{i-1}\right]B_{n+i},
\end{equation}
which reveals that \eqref{e:Kaneko} is simply the case $m=n+1$ of Lemma \ref{L:Kaneko}.

We can generalize Kaneko's identity in the following way:

\begin{theorem}
\begin{align}
\frac1{n+1}\sum_{i=0}^{n+1}2^{n+1-i}\binom{n+1}i \B_{n+i}&=(-1)^n,\notag\\
\frac1{n+1}\sum_{i=0}^{n+1}3^{n+1-i}\binom{n+1}i \B_{n+i}&=(-2)^{n-1}(n-4),\notag\\
\frac1{n+1}\sum_{i=0}^{n+1}4^{n+1-i}\binom{n+1}i \B_{n+i}
    &=(-1)^n\bigl(4^n+(2-\tfrac 43 n)3^n\bigr),\notag
\end{align}
and in general,
\begin{equation}
\frac1{n+1}\sum_{i=0}^{n+1}k^{n+1-i}\binom{n+1}i \B_{n+i}
   =\sum_{i=1}^{k-1}\bigl((2n+1)i-(n+1)k\bigr)i^n(i-k)^{n-1}.
 \label{e:genK}
\end{equation}
\end{theorem}
\begin{proof}
Using \eqref{e:Kkey}, we see that the left side of \eqref{e:genK}
is 
\[(B+k)^{n+1}B^n+B^{n+1}(B+k)^n.\] 
Setting $f(x)=x^m(x-k)^n$ in \eqref{e:b4}, we have
\[(B+k)^mB^n-(-1)^{m+n}B^m(B+k)^n=\sum_{i=0}^{k-1}\bigl((m+n)i-km\bigr)i^{m-1}(i-k)^{n-1}.\]
Setting $m=n+1$ gives \eqref{e:genK}.

\end{proof}

\bigskip
There are several interesting identities for Bernoulli numbers that actually hold for any 
two sequences
$(c_n)$ and $(d_n)$ related umbrally by $d^n=(c+1)^n$; i.e.,
\begin{equation}
\label{e:cd}
d_n=\sum_{i=0}^n\binom ni c_i.
\end{equation}
We note that \eqref{e:cd} may
inverted to give $c^n=(d-1)^n$; i.e.,
\[c_n=\sum_{i=0}^n(-1)^{n-i}\binom ni d_i.\]  By Theorem \ref{T:b3} (i),
\eqref{e:cd} holds with $c_n=B_n$, $d_n=(-1)^nB_n$.  We shall next
describe several pairs of sequences satisfying 
\eqref{e:cd}, and then give some identities for such sequences, which seem
to be new.

Since
\eqref{e:cd} is equivalent to
\[\egf{d_n}=e^x\egf{c_n},\] it is easy to find sequences satisfying
\eqref{e:cd} with simple exponential generating functions, though not all of our examples are of this form.

The derangement numbers $D_n$ satisfy $n!=\sum_{i=0}^n \binom ni D_i$ so \eqref{e:cd} holds
with $c_n=D_n$,
$d_n=n!$.

For any fixed nonnegative integer $m$, the Stirling numbers of the second kind $S(m,n)$ satisfy 
$n^m=\sum_{i=0}^n \binom ni i!\, S(m,i),$ so \eqref{e:cd} holds with $c_n=n!\,S(m,n)$, 
$d_n=n^m$.

The Euler numbers $E_n$ are defined by
$\egft{E_n}=\sech x$.
Let us define the ``signed tangent numbers" $T_n$ by $\tanh x=\egft {T_n}$.
Then since
$e^x\sech x=1+\tanh x$, we have that
\eqref{e:cd} holds with $c_n=E_n$, $d_n=\delta_n+T_n$.

The Genocchi numbers $g_n$ are defined by $\egft{g_n}=2x/(e^x+1)$. Then 
\[\frac{2xe^x}{e^x+1}=2x-\frac{2x}{e^x+1},\]
so \eqref{e:cd} holds with $c_n=g_n$, $d_n=2\delta_{n-1}-g_n$ (so that $d_1=g_1=1$).

The Eulerian polynomials $A_n(t)$ satisfy 
\[\egf {A_n(t)}=\frac{1-t}{1-te^{(1-t)x}}.\]
Then $A_0=1$, and $A_n(t)$ is divisible by $t$ for $n>1$. Let us set $\tilde 
A_n(t)=t^{-1}A_n(t)$ for $n>0$,
with $\tilde A_0(t)=1$. It is easy to check that 
\[e^{(1-t)x}\egf{\tilde A_n(t)}=\egf{A_n(t)},\]
so \eqref{e:cd} holds with $c_n=A_n(t)/(1-t)^n$, $d_n=\tilde A_n(t)/(1-t)^n$.

The Fibonacci numbers $F_n$ are defined by $F_0=1$, $F_1=1$, and $F_n=F_{n-1}+F_{n-2}$ for all 
integers $n$. It
is easily verified that for every fixed integer $m$, \eqref{e:cd} holds with $c_n=F_{m+n}$, 
$d_n=F_{m+2n}$, and also with
$c_n=F_{m-n}$,
$d_n=F_{m+n}$.

By the Chu-Vandermonde theorem, \eqref{e:cd} holds with \[c_n=(-1)^n 
\frac{(\alpha)_n}{(\beta)_n},\quad
d_n=\frac{(\beta-\alpha)_n}{(\beta)_n},\]
where $(\alpha)_n=\alpha(\alpha+1)\cdots(\alpha+n-1)$.

As Zagier observed \cite{kaneko}, it is easy to characterize the pairs of sequences satisfying 
\eqref{e:cd} with
$d_n=(-1)^nc_n$, which, as we shall see,  give analogues of Kaneko's
identity.  The condition, with  $c(x)=e^{cx}$ and $d(x)=e^{dx}$, is $e^x
c(x)=c(-x)$, which is equivalent to $e^{x/2}c(x)=e^{-x/2}c(-x)$;
i.e., $e^{x/2}c(x)$ is even. Thus it is easy to construct such sequences,   but not many seem
natural. In addition to the Bernoulli numbers, we have an example with the Genocchi  numbers $g_n$,
\[c(x)=\frac{2e^{-x/2}}{e^{x/2}+e^{-x/2}}=\frac 2{e^x+1}=\egf {\frac {g_{n+1}}{n+1}},\] 
and one with the Lucas numbers, $c_n=(-2)^{-n}(L_n+L_{2n})$, where $L_n=F_{n+1}+F_{n-1}$,

We now discuss the identities  which are consequences of \eqref{e:cd}.  Our first
identity generalizes Lemma  \ref{L:Kaneko}.

\begin{theorem}
\label{T:cd0}
Suppose that the sequences $c_n$ and $d_n$ satisfy \eqref{e:cd}. Then for all nonnegative 
integers $m$ and $n$,
\begin{equation}
\sum_{i=0}^m\binom mi c_{n+i}=\sum_{j=0}^n \binom nj(-1)^{n-j}d_{m+j}.
\label{e:cd0}
\end{equation}
\end{theorem} 

\begin{proof}
Let $c$ and $d$ be umbrae with $c^n=c_n$ and $d^n=d_n$. Then \eqref{e:cd} implies that 
$(c+1)^n=d^n$, so 
for any polynomial $f(x)$, we have $f(c+1)=f(d)$. Taking $f(x)=x^m(x-1)^n$ yields the
theorem.
\end{proof}

An  application of Theorem \ref{T:cd0} yields an interesting recurrence for 
Genocchi numbers. Let
$c_n$ be the Genocchi number $g_n$, so that, as noted above, $d_n=2\delta_{n-1}-g_n$. Then 
taking $m=n$ in
Theorem
\ref{T:cd0}, we have for $n>1$, 
\[\sum_{i=0}^n\binom ni g_{n+i}=-\sum_{i=0}^n \binom ni(-1)^{n-i}g_{n+i},\]
so
\[\sum_{i=0}^n(1+(-1)^{n-i})\binom ni g_{n+i}=0.\]
The only nonzero terms in the sum are those with $n-i$ even, so we may set $2j=n-i$ and divide 
by 2 to get the
recurrence
\begin{equation}
\label{e:seidel}
\sum_{j=0}^{\lfloor n/2\rfloor}\binom n{2j}g_{2n-2j}=0,\quad n>1.
\end{equation}
Equation  \eqref{e:seidel} is known as Seidel's recurrence (see, e.g., Viennot \cite{viennot}). 
It implies that $g_{2n}$ is
an integer, which is not obvious from the generating function (it is easily shown that 
$g_{2i+1}=0$
 for $i>0$), and it can also be used to derive a combinatorial interpretation for the Genocchi 
numbers. The reader can
check that the Genocchi analogue of Kaneko's identity alluded to 
before  Theorem \ref{T:cd0} is also
Seidel's recurrence in the form $\sum_{i=0}^n \binom ni g_{n+i}=0$ (in this form true for all 
$n$).

We now derive some further identities 
for sequences satisfying 
\eqref{e:cd}, of which the first generalizes Theorem \ref{T:cd0}.
\begin{theorem}
Suppose that the sequences $(c_n)$ and $(d_n)$ satisfy 
$\sum_{i=0}^n \binom ni c_i=d_n$.
Then for all $a$ and $b$, 
\begin{gather}
\sum_{i=0}^n\binom ai\binom b{n-i}d_i=\sum_{j=0}^n\binom aj\binom{a+b-j}{n-j}c_j
\label{e:cd1}
\\
\sum_{i=0}^n\binom ai\binom{2a-2i}{n-i}(-2)^i d_i
   =\sum_{j=0}^{\lfloor n/2\rfloor}\binom a{n-j}\binom{n-j}j(-2)^{n-2j}c_{n-2j}
\label{e:cd2}\\
\sum_{i=0}^n\binom ai\binom{2a-2i}{n-i}(-4)^i d_i 
   =(-1)^n\sum_{j=0}^n\binom aj\binom{2a-2j}{n-j}4^j c_j.
\label{e:cd3}
\end{gather}
\end{theorem}

\begin{proof}
Let $c$ be an umbra with $c^n=c_n$. Then \eqref{e:cd1}--\eqref{e:cd3} follow by substituting $c$ 
for $u$ in the 
following polynomial identities, where $v=1+u$:

\begin{gather}
\label{e:uv1}
\sum_{i=0}^n\binom ai\binom b{n-i}v^i=\sum_{j=0}^n\binom aj\binom{a+b-j}{n-j}u^j\\
\label{e:uv2}
\sum_{i=0}^n\binom ai\binom{2a-2i}{n-i}(-2v)^i
  =\sum_{j=0}^{\lfloor n/2\rfloor}\binom a{n-j}\binom{n-j}j(-2u)^{n-2j}\\
\label{e:uv3}
\sum_{i=0}^n\binom ai\binom{2a-2i}{n-i}(-4v)^i =
  (-1)^n\sum_{j=0}^n\binom aj\binom{2a-2j}{n-j}(4u)^j.
\end{gather}

We prove \eqref{e:uv1} by extracting the coefficient of $x^n$ in $(1+x)^b(1+vx)^a$ in two ways.
It is clear that this coefficient is given by the left side of \eqref{e:uv1}. But we also have
\[ (1+x)^b(1+vx)^a=(1+x)^b(1+x+ux)^a=(1+x)^{a+b}\left(1+\frac{ux}{1+x}\right)^a,\]
in which the coefficient of $x^n$ is easily seen to be given by the right side of  
\eqref{e:uv1}.

For \eqref{e:uv2}, we extract the coefficient of $x^n$ in 
\[\bigl((1+x)^2-2xv\bigr)^a=\bigl((1+x)^2-2x-2xu\bigr)^a
  =(1+x^2-2xu)^a.\]
For the left side we have
\begin{align*}
  \bigl((1+x)^2-2xv\bigr)^a&=(1+x)^{2a}\left(1-\frac{2xv}{(1+x)^2}\right)^a\\
  &=(1+x)^{2a}\sum_i\binom ai \frac{(-2xv)^i}{(1+x)^{2i}}\\
  &=\sum_i \binom ai x^i(1+x)^{2a-2i}(-2v)^i
\end{align*} 
and the coefficient of $x^n$ is the left side of \eqref{e:uv2}.

For the right side we have 
\begin{align*}
(1+x^2-2xu)^a&=\sum_i \binom ai (x^2-2xu)^i\\
  &=\sum_{i,j} \binom ai \binom ij x^{2j}(-2xu)^{i-j}\\
  &=\sum_{i,j} \binom ai \binom ij x^{i+j}(-2u)^{i-j}.\\
\end{align*}
Setting $i=n-j$ gives the right side of \eqref{e:uv2} as the coefficient of $x^n$.

For \eqref{e:uv3} we start with the identity   
\[\bigl((1+x)^2-4xv\bigr)^a
  =\bigl((1+x)^2-4x-4xu\bigr)^a=\bigl((1-x)^2-4xu\bigr)^a.\]
The coefficient of $x^n$ may be extracted from both sides as on the
left side of \eqref{e:uv2}.
\end{proof}

We note that  \eqref{e:uv1} is equivalent to a $_2F_1$ linear transformation and
\eqref{e:uv2} to a
$_2F_1$ quadratic transformation. Equation 
\eqref{e:uv3} is actually  a special case of
\eqref{e:uv1}; it can be obtained from \eqref{e:uv1} by replacing $a$ with $2a-n$ and $b$ with 
$n-a-\half$,
and simplifying.

The special case $a=-1$ of \eqref{e:cd1} is worth noting. It may be
written
\begin{equation}
\label{e:a=-1}
\sum_{i=0}^n \binom{b}{n-i} (-1)^i d_i =(-1)^n\sum_{j=0}^n
\binom{n-b}{n-j} c_j.
\end{equation}
If  we replace $n$ by $m+n$ in \eqref{e:a=-1} and then set $b=n$, it reduces to  \eqref{e:cd0}.

Next we prove a remarkable identity of Zagier \cite{zagier} for
Bernoulli numbers. Our proof is essentially an umbral version of
Zagier's. The reader may find it instructive to compare the two
presentations.  \def\BB{B^*}

\begin{theorem}
Let \[\BB_n=\sum_{r=0}^n \binom{n+r}{2r}\frac{B_r}{n+r}\] for
$n>0$. Then the value of $\BB_n$ for $n$ odd is periodic and is given
by
\smallskip
$$\offinterlineskip\vbox{\halign{\hfil#\hfil&&\vrule\hfil \hbox to 40pt{\hfil $#$\hfil}\cr
\vrule height 0pt width 0pt depth 5pt  $n$ {\rm(mod 12)}\kern 3pt&1&3&5&7&9&11\cr
\noalign{\hrule}
\vrule height 12pt width 0pt depth 0pt 
$B_n^*$&\tfrac34&-\tfrac14&-\tfrac14&\tfrac14&\tfrac14&-\tfrac34\cr
}}$$
\end{theorem}

\begin{proof}
Since $\binom{n+r}{2r}=\binom{n+r-1}{2r-1}\frac{n+r}{2r}$ for $r>0$, we have
\begin{align*}
2\sum_{n=1}^\infty \BB_n x^n
  &=2\sum_{n=1}^\infty \biggl[\frac 1n 
+\sum_{r=1}^n\binom{n+r-1}{2r-1}\frac{B^r}{2r}\biggr]x^n\\  
  &=-\log(1-x)^2-\log\left(1-B\frac x{(1-x)^2}\right) \\
  &= -\log\left((1-x)^2-Bx\right).
\end{align*}
Now let $g(u)=-\log(1-ux+x^2)$, so that $2\sum_{n=1}^\infty \BB_n x^n=g(B+2)$. Note that
\[g'(u)=\frac x{1-ux+x^2}.\]
Taking $k=4$ and $f(u)=g(u-2)$ in \eqref{e:b4}, we have
\begin{align*}
g(B+2)-g(-B-2)&=g'(-1)+g'(0)+g'(1)\\
  &=\frac x{1+x+x^2}+\frac x{1+x^2}+\frac x{1-x+x^2}\\[6pt]
  &=\frac{3x-x^3-x^5+x^7+x^9-3x^{11}}{1-x^{12}}.
\end{align*}
But 
$g(-B-2)=-\log\left((1+x)^2+Bx\right)=2\sum_{n=1}^\infty \BB_n(-x)^n$,
so
\[g(B+2)-g(-B-2)=4\sum_{n\text{ odd}}\BB_n x^n,\]
and the result follows.
\end{proof}

\section{Kummer Congruences}
\label{s:kummercong}
We say that a sequence $(u_n)$ of integers satisfies \emph{Kummer's congruence} for the
prime $p$ if for every integer $n$ and every $j\ge n$,
\begin{equation}\sum_{i=0}^n (-1)^{n-i}\binom n i u_{i(p-1)+j}\equiv
0\pmod{p^n}.\label{eq:K1}\end{equation} There are many  variations and generalizations of this
congruence, and we refer the reader to
\cite{geneul}, on which  most of this section is based, for more information and further
references.

If we set $j=n+k$, then \eqref{eq:K1} may be written umbrally as
\begin{equation}
(u^p-u)^nu^k\equiv 0 \pmod{p^n}
\label{e:Kummeru}
\end{equation}
for all $n,k\ge0$,
where $u^m=u_m$.

The result that we prove here shows that if  a sequence satisfies Kummer's congruence, then so 
does the
coefficient sequence of the reciprocal of its exponential generating function. Similar results 
apply to
products.

\begin{theorem}\label{T:Kummer} Let $(u_n)$ and $(v_n)$ be sequences of integers
satisfying 
\begin{equation}\left(\sumz n u_n \dpow xn\right)\left(\sumz n v_n \dpow xn\right)=1.
\label {eq:inv}\end{equation}
Then if $(u_n)$ satisfies Kummer's congruence for the prime $p$, so does $(v_n)$.
\end{theorem}
\begin{proof}
The relation \eqref{eq:inv} may be written umbrally as 
$(u+v)^n=0$, for $n>0$, where $u$ and $v$ are independent umbrae satisfying $u^m=u_m$ and
$v^n=v_n$, and this implies that if $f(x)$ is any polynomial with no constant term, then
$f(u+v)=0$. We shall prove by induction that if \eqref{e:Kummeru} holds for all
$n,k\ge0$ then 
\begin{equation}
(v^p-v)^n v^k\equiv 0\pmod{p^n}
\label{e:Kummerv}
\end{equation}
for all $n,k\ge0$.

The case $n=0$ of \eqref{e:Kummerv} is trivial. Now let 
$N$ be a positive integer and   $K$ a nonnegative integer, and suppose that 
\eqref{e:Kummerv} holds whenever $n<N$ and also when $n=N$ but $k<K$.
Thus 
\begin{align*}
0&=[(u+v)^p-(u+v)]^N (u+v)^K\\
&=[(u^p-u)+(v^p-v)+pR(u,v)]^N(u+v)^K\\
\intertext{for some polynomial $R(u,v)$ with integer coefficients,}
&=(v^p-v)^Nv^K+\text{other terms}.
\end{align*}
Here each other term is an integer times
$(u^p-u)^a(v^p-v)^b\bigl(pR(u,v)\bigr)^cu^dv^e$, where $a+b+c=N$, $d+e=K$, and either
$b<N$ or $b=N$, $c=0$, and $e<K$. Thus  by the inductive hypothesis and 
\eqref{e:Kummeru}, each of the
other terms is divisible by $p^{a+b+c}=p^N$, and therefore $(v^p-v)^Nv^K$ is also.
\end{proof}

As an example, we apply Theorem \ref{T:Kummer} to generalized Euler numbers. Recall that the
Euler numbers $E_n$ are defined by $\sech x=\sumz n E_n x^n/n!$ (so $E_n=0$ when $n$ is odd). We
define the generalized Euler numbers $\e mn$ by 
\[\sumz n \e mn \dpowp x{mn}=\left(\sumz n \dpowp x{mn}\right)^{-1},\]
so that $\e 2n=E_{2n}$. 
\begin{theorem}
\label{T:geneul}
Let $p$ be a prime and let $m$ be a positive integer such that $d=(p-1)/m$ is an
integer. Then for
$j\ge n/m$,
\[\sum_{i=0}^n (-1)^{n-i}\binom ni \e m{id+j}\equiv 0\pmod {p^n}.\]
\end{theorem}
\begin{proof}
Let us take $u_n=1$ and $v_n=\e m{n/m}$ if $m$ divides $n$, with $u_n=v_n=0$ otherwise.
Then the sequences $(u_n)$ and $(v_n)$ satisfy \eqref{eq:inv}, and $(u_n)$ satisfy
Kummer's congruence for $p$. Therefore $(v_n)$ does also.
\end{proof}

For example, if we take $m=4$ and $p=5$ in Theorem \ref{T:geneul}, then $d=1$ and we
have the congruence
\begin{equation}
\label{e:K5}
\sum_{i=0}^n (-1)^{n-i}\binom ni \e4{i+j}\equiv 0\pmod {5^n}
\end{equation}
for $j\ge n/4$.

By the same kind of reasoning we can prove a variation of Theorem \ref{T:Kummer} \cite{geneul}:

\begin{theorem}
\label{T:Kummer2}
Let the sequences $(u_n)$ and $(v_n)$ be related  by \eqref{eq:inv}, and suppose that for some 
integer $a$, 
\[\sum_{i=0}^n a^{n-i}\binom ni u_{ip+j}\equiv 0  \pmod {p^n}\]
for all nonnegative integers $j$ and $n$.
Then
\[\sum_{i=0}^n (-a)^{n-i}\binom ni v_{ip+j}\equiv 0  \pmod {p^n}.\qed\]
\end{theorem}

Next we prove a Kummer congruence for Bernoulli numbers.
A similar, but weaker, congruence was proved by Carlitz \cite{carlitzk} using a different method.

We call a
rational number \emph{$2$-integral} if its denominator is odd. If $a$ and $b$ are rational
numbers, then by $a\equiv b\pmod{2^r}$ we mean that $(a-b)/2^r$ is 2-integral. For example, 
$\half\equiv \tfrac52\pmod 2$.
We define $\rho_2(a)$ to be the largest integer for which $a/2^{\rho_2(a)}$ is
2-integral; so $\rho_2(\half)=-1$ and $\rho_2(\tfrac43)=2$.

In the proof of the next theorem we will use the fact that $2B_n$ is 2-integral
for all $n$, and that if $n$ is even and positive then $B_n\equiv\half\pmod 1$;
this follows easily by induction from the case $k=2$ of Theorem \ref{T:b3}
(iii).

\begin{theorem}
\label{T:KB}
For nonnegative integers $n$ and $j$,
\begin{equation}
\sum_{i=0}^n (-1)^{n-i}\binom ni B_{2i+2j}
\equiv 0\pmod{2^{\tau_{j,n}}},
\label{e:KB}
\end{equation} where
$\tau_{0,0}=0$, $\tau_{j,0}=\tau_{0,n}=-1$ for $j>0$ and $n>0$,
$\tau_{j,1}=1$ for $j\ge2$, and 
\[\tau_{j,n}=\min\left(2j-2, 2\left\lfloor\frac{3n-1}{2}\right\rfloor\right)\] for $n\ge2$
and
$j\ge1$. Moreover, the exponent in
\eqref{e:KB} is best possible if and only if $j\ne \lfloor (3n+1)/2\rfloor$. 
\end{theorem}

\begin{proof}
For simplicity we prove only the most interesting case, in which $n\ge2$ and $j\ge1$.
Let $B$ be the Bernoulli umbra, $B^n=B_n$, so the sum in \eqref{e:KB} is 
$(B^2-1)^n B^{2j}$.

Applying \eqref{e:bx3} with $k=2$ and $f(u)=(u^2-1)^{n}u^{2j}$, we obtain
\begin{equation}
(B^2-1)^n B^{2j}=2^{2j-1}(4B^2-1)^{n}B^{2j}+2^{2n-1}B^n(B+1)^n (2B+1)^{2j}.
\label{e:KB1}
\end{equation}

The first term on the right side of \eqref{e:KB1} is
$(-1)^n2^{2j-1}(B_{2j}-4nB_{2j+2}+\cdots)$.
Since $j>0$, this is congruent to
$(-1)^n 2^{2j-2} \pmod {2^{2j}}$ and thus $\rho_2(2^{2j-1}(4B^2-1)^{n}B^{2j})=2j-2$.

Next, let $g(u)=u^n(u+1)^n (2u+1)^{2j}$. To determine  $\rho_2\bigl (g(B)\bigr)$, we apply
\eqref{e:b2} in the form $g(B)=g(B-1)-g'(-1)$  and we find that
(since $n>1$)
\[g(B)=B^n(B+1)^n(2B+1)^{2j}=B^n(B-1)^n(2B-1)^{2j}.\]
We now apply \eqref{e:double} to $f(u)=u^n(u-1)^n(2u-1)^{2i}$ and we obtain
(since $n>1$)
\[g(B)=f(B)=f(2B)=2^nB^n(2B-1)^n(4B-1)^{2j}=(-2)^n\bigl(B_n-2nB_{n+1}+2K\bigr),\]
where $K$ is 2-integral. Thus if $n$ is even,
\[g(B)\equiv 2^nB_n \equiv 2^{n-1}\pmod {2^n},\]
 and if $n$ is odd 
\[g(B)\equiv 2^{n+1}nB_{n+1}\equiv 2^n\pmod {2^{n+1}}.\]
Thus $\rho_2(g(B))$ is $n-1$ if $n$ is even
and $n$ if $n$ is odd;  so in either case we have $\rho_2(g(B))=2\lfloor(n-1)/2\rfloor+1$. Thus the
power of 2 dividing the second term on the right side of \eqref{e:KB1} is
\[\rho_2(g(B))= (2n-1)+2\lfloor(n-1)/2\rfloor+1 
=2\lfloor(3n-1)/2\rfloor,\] 
and the congruence \eqref{e:KB} follows. It is clear that the exponent in \eqref{e:KB} 
is best possible if and
only if  $2j-2\ne 2\lfloor(3n-1)/2\rfloor$ and this is equivalent to the stated condition.
\end{proof}

We can use Theorem \ref{T:KB} to obtain congruences of a different kind for the Bernoulli numbers.
As noted earlier, for sequences $(c_n)$ and $(d_n)$ we have
$c_n=\sum_{i=0}^n(-1)^{n-i}\binom ni d_i$ if and only if $d_n=\sum_{i=0}^n \binom ni c_i$.
Let us fix $j>0$ and take $d_n=B_{2n+2j}$, so that $c_n=(B^2-1)^n B^{2j}$. Then we have
\[B_{2n+2j}=\sum_{i=0}^n\binom{n}{i}c_i.\] Moreover, it follows from Theorem 
\ref{T:KB} that if $i\ge(2j-1)/3$ and $i\ge2$ then $c_i\equiv 0\pmod
{2^{2j-2}}$, so we obtain the congruence
\begin{equation}
B_{2n+2j}\equiv\sum_{i=0}^M \binom ni c_i\pmod{2^{2j-2}},\label{e:Bcong}
\end{equation}
where $M=\max(\lfloor2(j-1)/3\rfloor,1)$. The cases $j=2,3,4$ of
\eqref{e:Bcong}, with simplifications obtained by reducing their coefficients, are
\begin{align*}
B_{2n+4}&\equiv-\frac1 {30}+\frac2{35}n\equiv \frac12+2n\pmod 4 \\
B_{2n+6}&\equiv \frac 1{42}-\frac2{35}n\equiv \frac{13}2+10n\pmod {16}\\
B_{2n+8}&\equiv -\frac 1{30}+\frac 6{55}n-\frac{2192}{5005}\binom{n}{2}
\equiv \frac{17}2+42n+ 48\binom n2 \pmod{64}.
\end{align*}
We note for use in the next
section simpler forms of the first two of these congruences: 

\begin{lemma}
\label{L:b2cong} Let $n$ be an even integer.
\begin{enumerate}
\item[(i)]
If $n\ge4$  then $B_n\equiv\frac12 +n\pmod 4$.
\openup 2\jot
\item[(ii)]
If $n\ge6$  then $B_n\equiv \frac12 +5n\pmod {16}$.\qed
\end{enumerate}
\end{lemma}
Of course, more direct proofs of this lemma are possible.
Similar congruences for Bernoulli numbers to other moduli have been given by Frame \cite{frame}.
Many congruences for generalized Euler numbers, obtained in this way from Kummer congruences,
can be found in \cite{geneul}.

\section{Median Genocchi numbers and Kummer congruences for Euler numbers}
\label{s:eulercong}
It follows from Theorem \ref{T:Kummer2} that the Euler numbers $E_n$  satisfy the congruence
\begin{equation*}
\sum_{i=0}^n\binom ni E_{2i+j}\equiv 0\pmod {2^n}.
\end{equation*}
However, Frobenius \cite{frobenius} (see also Carlitz \cite{carlitzk})
proved a much stronger congruence: the power of 2 dividing
\[\sum_{i=0}^n (-1)^{n-i}\binom ni E_{2i+j}\]
(for $j$ even) is the same as the power of 2 dividing $2^n n!$. Using the same approach as
Frobenius and Carlitz, we prove in this section a Kummer congruence for the numbers
$F_n={nE_{n-1}}$ in which the modulus is $2^{3n}$ or $2^{3n-1}$, and derive from it 
Frobenius's congruence.

A special case of our result gives a divisibility property for the \emph{median Genocchi numbers}
(also 
called \emph{Genocchi numbers of the second kind}). These numbers may be defined by 
\begin{equation*}
H_{2n+1}=\sum_{k=0}^n \binom nk g_{n+k+1},
\end{equation*}
where the $g_i$ are the Genocchi numbers, defined by $\sumz i g_i x^i/i!=2x/(e^x+1)$.
(In combinatorial investigations, 
the notation $H_{2n+1}$ is usually used for what in our notation is $|H_{2n+1}|= (-1)^n
H_{2n+1}$.) The connection between  median Genocchi numbers and the numbers
$F_n=nE_{n-1}$ is given by  the following result, due to Dumont and Zeng
\cite{dz}.

\begin{lemma}
\label{L:median}
\[2^{2n}H_{2n+1}=\sum_{i=0}^n(-1)^{n-i}\binom ni F_{2i+1}=
  \sum_{i=0}^n (-1)^{n-i}\binom ni (2i+1)
E_{2i}.\]
\end{lemma}
\begin{proof}
Let us define the umbrae $g$ and $F$ by $g^n=g_n$ and $F^n=F_n=nE_{n-1}$,
so that $e^{Fx}=x\sech x$.
Then
\[e^{gx}=\frac{2x}{e^x+1}=\frac{2xe^{-x/2}}{e^{x/2}+e^{-x/2}}=2e^{-x/2}\cdot\frac x2 \sech \frac x2 =2e^{-x/2}e^{\frac F2 x}
=2e^{\frac 12 (F-1)x}.\]
Thus $g^n=2\left(\frac{F-1}2\right)^n$,
and it follows by linearity that 
\begin{align*}
2^{2n}H_{2n+1}&= 2^{2n}(g+1)^n
g^{n+1}=2^{2n}\cdot 2\left(\frac{F+1}2\right)^n\left(\frac{F-1}2\right)^{n+1}
\\
  &=(F-1)(F^2-1)^n=F(F^2-1)^n,
\end{align*}
since $F^m=0$ for $m$ even.
\end{proof}

Barsky \cite{barsky} proved a conjecture of Dumont that $H_{2n+1}$ is divisible by $2^{n-1}$. More
precisely, Barsky proved that for $n\ge3$, $H_{2n+1}/2^{n-1}$ is congruent to 2 modulo 4 if $n$ is
odd and is congruent to 3 modulo 4 if $n$ is even. Kreweras \cite{kreweras} gave a combinatorial
proof of Dumont's conjecture, using a combinatorial interpretation of $H_{2n+1}$ due to Dumont
\cite[Corollaire 2.4]{dr}. A $q$-analogue of Barsky's result was given by Han and Zeng \cite{hz}.

It is interesting to note that (as pointed out in \cite{sloane}), a combinatorial interpretation of
the numbers $H_{2n+1}/2^{n-1}$ was given in 1900 by H. Dellac \cite{dellac}. Dellac's 
interpretation may be described as follows: 
We start  with a $2n$ by $n$ array of cells and  consider the set $D$  of cells in
rows $i$ through $i+n$ of column $i$, for $i$ from 1 to $n$. Then $H_{2n+3}/2^n$ is the number of
subsets of $D$ containing two cells in each column and one cell in each row. Dellac did not give
any formula for these numbers, but he did compute them for $n$ from 1 to 8. Dellac's
interpretation can be derived without too much difficulty from Dumont's combinatorial
interpretation, but it is not at all clear how Dellac computed these numbers.

\begin{theorem}\label{T:xsech}
Let $F_n=nE_{n-1}$, so that $\egft {F_n}=x\sech x$. Let $2^{\mu_{j,n}}$ be the highest power of 
$2$ dividing
\[\sum_{i=0}^n (-1)^{n-i}\binom ni F_{2i+j},\]
where $j$ is odd.
Then $\mu_{j,0}=0$, $\mu_{j,1}=2$, and for $n>1$, $\mu_{j,n}=3n$ if $n$ is odd and 
$\mu_{j,n}=3n-1$ if $n$ is even.
\end{theorem}

\begin{proof}
We have
\[
\sech x - 2e^{-x}=\frac{2}{ e^{x}+e^{-x}}-\frac{2+2e^{-2x}}{e^{x}+e^{-x}}
=-\frac{2e^{-2x}}{e^{x}+e^{-x}}=-\frac{2e^{-x}}{e^{2x}+1}
=-2\frac{e^x-e^{-x}}{e^{4x}-1}.
\]
Therefore
\begin{equation}x\sech x= 2xe^{-x}-\sinh x\cdot B(4x)
\label{e:sech}
\end{equation}
where $B(x)=x/(e^x-1)$ is the Bernoulli number generating function.

Now let us define the umbrae $F$, $A$, $B$, and  $C$ by 
\begin{align*}
e^{Fx}&=x\sech x\\
e^{Ax}&=xe^{-x}=\sum_{n=1}^\infty (-1)^{n-1}n\dpow xn\\
e^{Bx}&=B(x)=\frac{x}{e^x-1}\\
e^{Cx}&=\sinh x=\sum_{n\ \text{odd}}\dpow xn.\\
\end{align*}
Then  from \eqref{e:sech} we have
\begin{equation}
\label{e:ABH}
F^n=2A^n-(4B+C)^n.
\end{equation}

We want to find the power of 2 dividing  $F^j(F^2-1)^n$. It follows from \eqref{e:ABH} that
\begin{equation}
F^j(F^2-1)^n=2A^j(A^2-1)^n-(4B+C)^j\bigl((4B+C)^2-1\bigr)^n.
\label{e:FABH}
\end{equation}
First note that for any polynomial $p$, $p(A)=p'(-1).$ Therefore,
\begin{equation}
2A^j(A^2-1)^n=\begin{cases} 2(-1)^{j-1}j&\text{if $n=0$}\\
  4(-1)^{j-1}&\text{if $n=1$}\\
  0&\text{if $n\ge 2$}\\
\end{cases}
\label{e:A}
\end{equation}
We note also that $C^j(C^2-1)^n=0$ for all integers $j\ge 0$ and $n\ge1$.
Then
\begin{align}(4B+C)^j\bigl((4B+C)^2-1\bigr)^n&=(4B+C)^j(16B^2+8BC+C^2-1)^n\notag\\
  &=2^{3n}(4B+C)^j(2B^2+BC)^n.\label{e:BH} 
\end{align}

We now need to determine the power of 2 dividing $(4B+C)^j(2B^2+BC)^n$. Since $F^j(F^2-1)^n=0$ 
if $j$ is even,
we may assume that $j$ is odd.
Since $2B_i$ is 2-integral, we have
\begin{equation*}
(4B+C)^j(2B^2+BC)^n\equiv C^j
(B^nC^n+2nB^{n+1}C^{n-1})\pmod2.
\end{equation*}
Using the facts that $B_i=0$ when $i$ is odd and greater than 1, and that 
\[C^i=\cases 1&\text{if $i$ is odd}\\ 0&\text{if $i$ is even,}\endcases\]
we find that if  $n=0$  then
\[(4B+C)^j(2B^2+BC)^n\equiv B_0\equiv1\pmod 2,\]
if  $n$ is even and positive then
\[(4B+C)^j(2B^2+BC)^n\equiv B_n\equiv\half\pmod 1,\]
and if  $n$ is odd then
\[(4B+C)^j(2B^2+BC)^n\equiv 2nB_{n+1}\equiv 1\pmod 2.\]
The theorem then follows from these congruences, together with 
\eqref{e:FABH}, \eqref{e:A}, and \eqref{e:BH}.
\end{proof}

By taking more terms in the expansion of $(4B+C)^j(2B^2+BC)^n$, we can get congruences modulo
higher powers of 2. 
For example, if  $n$ is even then 
\[(4B+C)^j(2B^2+BC)^n\equiv B_n+4\binom n2 B_{n+2}\pmod 8.\]
Applying Lemma \ref{L:b2cong} (ii),  we find that if $n\ge 6$ then 
\[B_n+4\binom n2 B_{n+2}\equiv \half+n^2(n+3)\pmod 8.\]
Since $n$  even implies $n^2\equiv 2n\pmod 8$ and $4n\equiv0 \pmod 8$, this simplifies to 
$\half-2n\pmod 8$.

Similarly, if $n$ is odd then 
\[(4B+C)^j(2B^2+BC)^n\equiv (2n+4j)B_{n+1}+4\binom{n}{3}\pmod 8.\]
If $n\ge 3$ then Lemma \ref{L:b2cong} (i) gives $(2n+4j)B_{n+1}\equiv (n+2j)(3+2n)\pmod 8$, and
we may easily verify that for $n$ odd,  $4\binom{n}{3}\equiv 2n-2 \pmod 8$.
Using the fact that $n$ odd implies  $n^2\equiv 1 \pmod 8$ and $4n\equiv 4
\pmod8$,  we obtain
\[(2n+4j)B_{n+1}+4\binom{n}{3}\equiv 4+2j+n \pmod 8.\]

Therefore we may conclude that (for $j$
odd) if
$n$ is even and $n\ge 6$ then
\[2^{-(3n-1)}\sum_{i=0}^n (-1)^{n-i}\binom ni F_{2i+j}\equiv 4n -1 \pmod {16},\]
and if $n$ is odd and $n\ge 3$ then 
\[2^{-3n}\sum_{i=0}^n (-1)^{n-i}\binom ni F_{2i+j}\equiv 4 - 2j -n \pmod 8.\]
In particular, we get a refinement of Barsky's theorem: 
If $n$ is even and $n\ge 6$ then $H_{2n+1}/2^{n-1}\equiv 4n-1\pmod {16}$, and if $n$ is odd and
$n\ge 3$ then $H_{2n+1}/2^n\equiv 2-n\pmod 8$.  It is clear that by the same method we could
extend these congruences to any power of 2.

Next we derive Frobenius's congruence  from
Theorem
\ref{T:xsech}. (This derivation is  similar
to part of Frobenius's original proof.) Define the umbra
$E$ by
$e^{Ex}=\sech x$, so $F^n=nE^{n-1}$. First note that if $j$ is even then $E_j=F_{j+1}/(j+1)$ is
odd, so Frobenius's congruence holds for $n=0$.
From  $F^n=nE^{n-1}$, it follows that for any polynomial
$p$, we have $p(F)=p'(E)$. Let us take  $p(u)=u^{j+1}(u^2-1)^n$, where $j$ is even and $n\ge1$.
Then 
$p'(u)=(j+1)u^{j}(u^2-1)^n+2nu^{j+2}(u^2-1)^{n-1}$.
By Theorem \ref{T:xsech}, we have $p(F)\equiv 0\pmod {2^{3n-1}}$, so
\[E^j(E^n-1)^n\equiv -\frac{2n}{j+1}E^{j+2}(E^2-1)^{n-1}\pmod{2^{3n-1}}.\]
By induction on $n$, the power of 2 dividing $E^{j+2}(E^2-1)^{n-1}$ is equal to the power of
2 dividing $2^{n-1}(n-1)!$, and  Frobenius's result follows.

In view of Theorem \ref{T:xsech}, it is natural to ask whether there are analogous congruences 
for generalized Euler
numbers. There seem to be many possibilities, but the most attractive is given by the following 
conjecture: 
Define numbers $\f mn$ by
\[\sumz n \f mn \frac {x^{(2n+1)m}}{\bigl((2n+1)m\bigr)!}=\dpow x m\bigg/\sumz n \dpowp 
x{2nm}.\] 
(Thus $\f 1n$ is $F_{2n+1}$ as defined above.) Let $2^{\mu_{j,n,t}}$ be the highest power of 2 
dividing
\[\sum_{i=0}^n (-1)^{n-i}\binom ni \f{2^t}{i+j}.\] Then for $t\ge1$, we have $\mu_{j,0,t}=0$,
$\mu_{j,1,t}=4$, and for $n >1$, $\mu_{j,n,t}=\lfloor \frac{7n}2 \rfloor-1$, except when $t=1$, 
$n\equiv 2\pmod 4$, and
$n\ge 6$. 

In the exceptional case, $\mu_{j,n,1}= \frac{7n}2 +2+\rho_2(j+\vartheta_n)$, where $\vartheta_n$ is
some integer or 2-adic  integer. The first few  values of $\vartheta_n$ (or
reasonably good 2-adic approximations to them) are $\vartheta_6=118$, $\vartheta_{10}=7$, 
$\vartheta_{14}=2$,
$\vartheta_{18}=13$, $\vartheta_{22}=32$, and $\vartheta_{26}=27$.

By way of illustration, $\mu_{0,1,t}=4$ for $t\ge1$ is equivalent to the (easily proved) 
assertion that $1+\binom{3\cdot
2^t}{2^t}$ is divisible by 16 but not by 32.

\section{Bell numbers}
\label{s:bell}
The Bell numbers $B_n$ are defined by the exponential
generating function
\begin{equation}
  B(x)= \egf {B_n}=e^{e^x-1}.\label{e:bell}
\end{equation}
Although we are using the same notation for the Bell numbers that we
used for Bernoulli numbers, there should be no confusion.  Rota \cite{partitions} proved
several interesting properties of the Bell numbers using umbral
calculus in his fundamental paper. Here we prove a well-known congruence of Touchard for Bell
numbers and a generalization due to Carlitz.

Differentiating \eqref{e:bell}
gives
$B'(x)=e^xB(x)$, so 
\[B_{n+1}=\sum_{k=0}^n \binom nk B_k.\]
With the Bell umbra $B$, given by $B^n=B_n$, this may be written
$B^{n+1}=(B+1)^n.$
 Then by linearity, for any polynomial $f(x)$ we have
\begin{equation}
Bf(B)=f(B+1)
\label{e:shift} 
\end{equation}

A consequence of \eqref{e:shift}, easily proved by induction, is that for any 
polynomial $f(x)$ and any nonnegative integer $n$, 
\begin{equation}
B(B-1)\cdots (B-n+1)f(B)=f(B+n).
\label{e:nshift}
\end{equation}
(A $q$-analogue of \eqref{e:nshift} has been given by Zeng [\ref{zeng}, Lemma 8].) If we take 
$f(x)=1$ in \eqref{e:nshift},
we obtain (since
$B_0=1$)
\begin{equation}
B(B-1)(B-2)\cdots(B-n+1)=1.
\label{e:rota}
\end{equation}
We note that Rota took \eqref{e:rota} as the \emph{definition}
of the Bell umbra and derived  \eqref{e:shift} and \eqref{e:bell} from it.

As an application of these formulas, we shall prove Touchard's congruence for the Bell numbers 
[\ref{touchard1},
\ref{touchard2}]. 

If $f(x)$ and $g(x)$ are two polynomials in $\Z[x]$, then by $f(x)\equiv g(x)\pmod p$, we mean
that $f(x)-g(x)\in p\Z[x]$. We first  recall two elementary facts about congruences for
polynomials modulo a prime $p$. First we have Lagrange's congruence,
$x(x-1)\cdots(x-p+1)\equiv x^p-x\pmod p$. Second, if
$g(x)\in \Z[x]$ then $g(x+p)-g(x)\equiv 0\pmod p$.

\begin{theorem}
For any prime $p$ and any nonnegative integer $n$, 
\[B_{n+p}-B_{n+1}-B_n\equiv 0\pmod p.\]
\end{theorem}

\begin{proof} By Lagrange's congruence, 
\begin{equation*}
 (B^p-B-1)B^n\equiv 
  \bigl(B(B-1)\cdots (B-p+1) -1\bigr)B^n\pmod p.
\end{equation*}
By \eqref{e:nshift},
\[\bigl(B(B-1)\cdots (B-p+1)-1\bigr)B^n=(B+p)^n-B^n.\]
Since  $(x+p)^n-x^n\equiv 0\pmod p$, $p$ divides $(B+p)^n-B^n$.
\end{proof}

Next we prove a generalization of Touchard's congruence analogous to a Kummer congruence,  due 
to
Carlitz \cite{carlitzbell}.

\begin{theorem}
\label{T:bellk}
For any prime $p$ and any nonnegative integers $n$ and $k$,
\[(B^p-B-1)^kB^n\equiv0 \pmod {p^{\ceil {k/2}}}.\]
\end{theorem}

\begin{proof}
Let $L(x)$ be the polynomial $x(x-1)\cdots(x-p+1) -1$. First we show that it suffices
to prove that for any  polynomial $f(x)\in \Z[x]$,
$L(B)^kf(B)\equiv 0 \pmod {p^{\ceil {k/2}}}$. To see this, note that we may write $L(x)=x^p-x-1 
- p R(x)$, where $R(x)\in
\Z[x]$. Then 
$(B^p-B-1)^kB^n=(L(B)+pR(B))^k B^n=\sum_{i=0}^k \binom ki p^i L(B)^{k-i}R(B)^i B^n$, and our 
hypothesis will show that 
$p^i L(B)^{k-i}R(B)^i B^n$ is divisible by $p$ to the power 
$i+\ceil{(k-i)/2}=\ceil{(k+i)/2}\ge\ceil{k/2}$.

We now prove by induction on $k$ that for any polynomial $f(x)\in \Z[x]$,
\[L(B)^kf(B)\equiv 0 \pmod {p^{\ceil {k/2}}}.\] The assertion is trivially true for $k=0$. 
For the induction step, note that we  may write $L(x+p)=L(x)+pJ(x)$, where $J(x)\in \Z[x]$, and 
recall that by 
\eqref{e:nshift}, $L(B)g(B)=g(B+p)-g(B)$ for any polynomial $g$.
Then for any $f(x)\in \Z[x]$ we have for $k>0$
\begin{align*}
L(B)^k f(B)&= L(B)\cdot L(B)^{k-1}f(B)=L(B+p)^{k-1}f(B+p)-L(B)^{k-1}f(B)\\
  &=\bigl(L(B)+pJ(B)\bigr)^{k-1}f(B+p)-L(B)^{k-1}f(B)\\
  &=L(B)^{k-1}\bigl(f(B+p)-f(B)\bigr)+\sum_{i=1}^{k-1}p^i\binom{k-1}i L(B)^{k-1-i}J(B)^i f(B+p).
\end{align*}
We show that each term of the last expression  is divisible by $p^{\ceil {k/2}}$. Since 
$f(x+p)-f(x)=ph(x)$ for some
$h(x)\in
\Z[x]$, 
we have $L(B)^{k-1}\bigl(f(B+p)-f(B)\bigr)=pL(B)^{k-1}h(B)$, which by induction is divisible by
$p$ to the power $1+\ceil{(k-1)/2}\ge\ceil{k/2}$. By induction also, the $i$th term in the sum 
is divisible by
$p$ to the power $i+\ceil{(k-1-i)/2}=\ceil{(k-1+i)/2}\ge\ceil{k/2}$. This completes the proof.
\end{proof}

Theorem \ref{T:bellk} can be extended in several ways (in particular, the modulus can be 
improved); see Lunnon, Pleasants, and  Stephens \cite{bellcong}
and Gessel \cite{comb}, which both use umbral methods (though the latter is primarily combinatorial).
Another congruence for Bell numbers, also proved umbrally, was given by Gertsch and Robert
\cite{gr}.

\bigskip
{\noindent\small\textsc{Acknowledgment:}
I would like to thank Jiang Zeng for his hospitality during my visit to the University of Lyon 1,
where most of this paper was written. I would also like to thank Timothy Chow, Hyesung Min, and an anonymous
referee for pointing out several errors in earlier versions of this paper.}

\end{document}